\documentclass[12pt, a4paper]{article}
\usepackage[latin1]{inputenc}
\usepackage[english]{babel}
\usepackage{indentfirst}
\usepackage{amsxtra}
\usepackage{amsmath}
\usepackage{amsthm}
\usepackage{amstext}
\usepackage{amsfonts}
\usepackage{textcomp}
\usepackage{amssymb}
\usepackage{amscd}
\usepackage[all]{xy}
\usepackage[dvips]{graphicx}
\usepackage{setspace}
\usepackage[all]{xy}

\newcommand {\demo}{\hskip -0.6cm{\bf Proof.  }}
\newcommand {\fim}{\hfill{$\square$}\vskip 1pc}

\newcommand {\R}{\mathbb{R}}
\newcommand {\N}{\mathbb{N}}

\newcommand {\X}{\mathbb{X}}

\newtheorem{teorema}{Theorem}[section]
\newtheorem{lema}[teorema]{Lemma}
\newtheorem{corolario}[teorema]{Corollary}
\newtheorem{definicao}[teorema]{Definition}
\newtheorem{proposicao}[teorema]{Proposition}
\newtheorem{exemplo}[teorema]{Example}
\newtheorem{obs}[teorema]{Remark}

\begin{document}

\onehalfspace

\title{Branching systems and representations of Cohn-Leavitt path algebras of separated graphs}

\maketitle
\begin{center}
%\vspace{5mm}
{\large Daniel Gonçalves\footnote{partially supported by CNPq} and Danilo Royer}\\
\end{center}
\vspace{8mm}

\begin{abstract}
 We construct for each separated graph $(E,C)$ a family of branching systems over a set $\X$ and show how each branching system induces a representation of the Cohn-Leavitt path algebra associated to $(E,C)$ as homomorphisms over the module of functions in $\X$. We also prove that the abelianized Cohn-Leavitt path algebra of a separated graph with no loops can be written as an amalgamated free product of abelianized Cohn-Leavitt algebras that can be faithfully represented via branching systems.
\end{abstract}

\section{Introduction}
\label{intro}
Cohn-Leavitt path algebras were introduced recently, see \cite{ara2}, \cite{ara1}, as generalizations of Leavitt path algebras, which in turn arised as algebraic analogues of graph C*-algebras (see \cite{Abrams}, \cite{Abrams2}). Both Leavitt path algebras and graph C*-algebras have been the focus of intense research in the last fiften years (see \cite{Abrams}, \cite{Abrams2}, \cite{PPM}, \cite{FLR}, \cite{GR1}, \cite{KPR} for a few examples) and nowadays the literature on then is extensive. This is not the case for Cohn-Leavitt path algebras, for which the greater generality of its definition allows for algebras that are very different from the usual Leavitt path algebras. In particular, in \cite{ara1} it is shown that any conical abelian monoid occurs as the monoid of isomorphism classes of
finitely generated projective modules over a Leavitt algebra of a separated graph (what, by the results of \cite{AMP}, is not true for the class of Leavitt algebras of non-separated graphs). Actually, although Cohn-Leavitt path algebras of separated graphs incorporate the usual Leavitt path algebras (for a particular separation), they behave quite differently from the later since the range projections associated to different edges need not commute. So, as expected, Leavitt path algebras and Cohn-Leavitt path algebras of separated graphs do not share all the same properties, and many results on  Leavitt path algebras still have to be extended. In particular, the results of \cite{rg1} regarding representations of Leavitt path algebras arising from branching systems still need Cohn-Leavitt path algebras versions.

Branching systems arise in many areas in mathematics, see \cite{uniteq}, and can be used to produce and study representations of Leavitt path algebras (see \cite{rg1}). More specifically, faithful representations of Leavitt path algebras can be obtained via branching systems and for certain graphs the study of representations, up to unitary equivalence, can be reduced to the study of representations arising from branching systems, see \cite{rg1}, \cite{uniteq}. It is our goal in this paper to extend some of the results concerning representations of Leavitt path algebras arising from branching systems (see \cite{rg1}) to Cohn-Leavitt path algebras.

We now describe precisely what we will do in this paper: In section 2 we define branching systems of separated graph and show how they induce representations of the associated Cohn-Leavitt path algebra in the algebra of homomorphisms over the module of functions over a set $\X$. Next, in section 3, we prove the existence of branching systems for any separated graph via a constructive argument and use the representations obtained to show a few properties of elements in the Cohn-Leavitt path algebra. We focus on the injectivety of the representations from branching systems in section 4. For this we must look into abelianized Cohn-Leavitt path algebras, since the image of the range projections of different edges under representations arising from branching system always commute, and we give some examples of separated graphs for which the representations of the abelianized algebra arising from the branching systems constructed in section 3 are injective. In particular, this includes graphs with no loops where all edges have the same source and the range map is injective. We then finalize the section showing that any Cohn-Leavitt path algebra can be written as an amalgamated free product of Cohn-Leavitt path algebras over graphs where all edges have the same source. Before we proceed we recall the definition of Cohn-Leavitt path algebras below.

A separated graph is a pair $(E,C)$ where $E=(E^0, E^1, r, s)$ is a directed graph, (that is, $E^0$ is a set of vertices, $E^1$ is a set of edges and $r,s:E^1\rightarrow E^0$ are the range and source maps), and $C=\bigcup\limits_{v\in E^0}C_v$, where each $C_v$ is a partition of $s^{-1}(v)$ into pairwise disjoint nonempty sets, for each non-sink $v$. Let $C_{fin}$ be the set of all finite sets $Y\in C$ and $\text{Path} (E)$ the set of all paths in $E$.

%We let $(E^1)^*$ denote the set of formal symbols $\{e^*:e\in E^1\}$ and for $\alpha:= e_1e_2 \ldots e_n \in E^n$ we define $\alpha^*:= e_n^*e_{n-1}^* \ldots e_1^*$ . We also define $v^*=v$ for all $v\in E^0$.

\begin{definicao}\label{defleviatt}(As in \cite{ara1}). Let $(E,C)$ be a separated graph, let $S\subseteq C_{fin}$ and, and $K$ be a field. The {\it Cohn-Leavitt algebra of the triple $(E,C,S)$}, denoted $L_K(E,C,S)$, is the universal $K$-algebra generated by a set $\{v:v\in E^0\}$, of pairwise orthogonal idempotents, together with a set $\{e,e^*:e\in E^1\}$ of elements satisfying:
\end{definicao}
\begin{itemize}
\item[(E1)] $s(e)e=er(e)=e$ for all $e\in E^1$
\item[(E2)] $r(e)e^*=e^*s(e)=e^*$ for all $e\in E^1$
\item[(SCK1)] $e^*f= \delta_{e,f} r(e)$ for all $e,f \in Y$, for each $Y\in C$
\item[(SCK2)] $v=\sum\limits_{e\in X}ee^*$ for every set $X\in S\cap C_v$, for each non-sink $v\in E^0$.
\end{itemize}

Following \cite{ara1}, $L_K(E,C,C_{fin})$ is the {\it Leavitt path algebra of the separated graph $(E,C)$}, and $L_K(E,C,\emptyset)$ is the Cohn path algebra of the graph $(E,C)$. If $C_v=s^{-1}(v)$ for each non-sink $v$ then $L_K(E,C,C_{fin})$ is the Leavitt path algebra of the directed graph $E$ (see \cite{Abrams}, \cite{Abrams2} for Leavitt path algebras of a graph.)

\section{$(E,C,S)$-algebraic branching systems}

In this section we will define $(E,C,S)$-algebraic branching systems associated to a triple $(E,C,S)$ and we will show how these $(E,C,S)$-algebraic branching systems induce representations of the associated Cohn-Leavitt path algebra, in the $K$ algebra of the homomorphisms in the module of functions over a set $\X$.

We start with the definition of an $(E,C,S)$-algebraic branching system:

\begin{definicao}\label{brancsystem}
Let $(E,C)$ be a separated graph and $S\subseteq C_{fin}$. Let $\X$ be a set and let $\{R_e\}_{e\in E^1}$, $\{D_v\}_{v\in E^0}$ be families of subsets of $\X$ such that:
\begin{enumerate}
\item $R_e\cap R_d= \emptyset$ for each $d,e\in Y$ with $d\neq e$, $Y\in C$,
\item $D_u\cap D_v= \emptyset$ for each $u,v\in E^0$ with $u\neq v$,
\item $R_e\subseteq D_{s(e)}$ for each $e\in E^1$.
\item $D_v=\bigcup\limits_{e\in Y}R_e$\,\,\,\,\, if\,\,\,\,\, $Y\in S\cap C_v$, for each non-sink $v\in E^0$.
\item for each $e\in E^1$, there exists a bijective map $f_e:D_{r(e)}\rightarrow R_e$.
\end{enumerate}

A set $\X$, with families of subsets $\{R_e\}_{e\in E^1}$, $\{D_v\}_{v\in E^0}$, and maps $f_e$ as above, is called an $(E,C,S)$- algebraic branching system, and we denote it by $(\X,\{R_e\}_{e\in E^1}, \{D_v\}_{v\in E^0}, \{f_e\}_{e\in E^1})$, or when no confusion arises, simply by $\X$.
\end{definicao}

Next, fix an $(E,C,S)$-algebraic branching system $\X$. Let $M$ be the $K$ module of all functions from $\X$ taking values in $K$ and let $Hom_K(M)$ denote the $K$ algebra of all homomorphisms from $M$ to $M$ (with multiplication given by composition of homomorphisms and the other operations given in the usual way).

Now, for each $e\in E^1$ and for each $v\in E^0$, we will define homomorphisms $S_e$, $S_e^*$ and $P_v$ in $Hom_K(M)$.

Let $S_e$ be defined as follows:
$$\left( S_e \phi \right) (x) =
\begin{cases} \phi(f_e^{-1}(x)),  \text { if }  x\in R_e \\ 0,  \text{ if } x\notin R_e \end{cases}, $$ where $\phi$ is a function in $M$.

In order to simplify notation, in what follows we will make a small abuse of the characteristic function symbol and denote the above homomorphism by:
$$S_e \phi=\chi_{R_e}\cdot \phi\circ f_e^{-1}.$$

In a similar fashion to what is done above, and making the same abuse of the characteristic function symbol, we define the homomorphism $S_e^*$ by $$S_e^* \phi=\chi_{D_{r(e)}}\cdot \phi\circ f_e,$$ where $\phi \in M$.

Finally, for each $v\in E^0$, and for $\phi \in M$, we define $P_v $ by $$P_v \phi=\chi_{D_v} \cdot \phi,$$
that is, $P_v$ is the multiplication operator by $\chi_{D_v}$, the characteristic function of $D_v$.

 \begin{teorema}\label{rep} Let $\X$ be an $(E,C,S)$- algebraic branching system. Then there exists a representation (that is, an algebra homomorphism) $\pi: L_K(E,C,S)\rightarrow Hom_K(M)$ such that
$$\pi(e)= S_e, \text{ } \pi(e^*)= S_e^* \text{ and } \pi(v)=P_v,$$ for each $e\in E^1$ and $v\in E^0$.
  \end{teorema}

  \demo As in \cite{rg1}

\fim

\begin{obs}\label{Mvanishing} Notice that theorem \ref{rep} still holds if we change the module $M$ of all functions from $\X$ to $K$ for the module of all functions from $\X$ to $K$ that vanish in all, but a finite number of points, of $\X$.
\end{obs}

In the next section we consider the question of existence of $(E,C,S)$-algebraic branching systems(and their induced representations) for any given graph $E$.

\section{Existence of $(E,C,S)$-algebraic branching systems}

Let $(E,C)$ be a separated graph, and $S\subseteq C_{fin}$. Next we show that there exists an $(E,C,S)$-algebraic branching system. Our proof is constructive and one can actually obtain a great number of $(E,C,S)$-algebraic branching systems following the ideas below.

\begin{teorema}\label{existencebrancsys}
Let $(E,C)$ be a separated graph, with $E^0$ and $E^1$ countable, and $S\subseteq C_{fin}$. Then there exists an $(E,C,S)$-branching system $\X$, with $\X\subseteq \R$, such that:
\begin{enumerate}
\item $R_e\cap R_f\neq \emptyset$ for each $e\in X$, $f\in Y$, and $X,Y\in C_v$ with $X\neq Y$.
\item for each $X\in C_v\setminus S$, $\bigcup\limits_{e\in X}R_e\subsetneq D_v$.
\item for $X,Y\in C_v\setminus S$ with $X\neq Y$ it holds that $\bigcup\limits_{e\in X}R_e\neq \bigcup\limits_{f\in Y}R_f$.

\end{enumerate}

\end{teorema}

\demo Since $E^0$ is countable then $E^0=\{v_i\}_{i=0}^N$, case if $E^0$ is finite, or $E^0=\{v_i\}_{i=0}^\infty$. For each $i$, let $D_{v_i}$ be the interval $[i,i+1)\subseteq\R$.

From now on, fix a non-sink $v\in E^0$. Since $E^1$ is finite or countable then  $C_v=\{Y_j\}_{j=1}^M$ or $C_v=\{Y_j\}_{j\in \N}$, and $Y_j=\{e_i^j\}_{i=1}^{K_j}$ or $Y_j=\{e_i^j\}_{i=1}^\infty$ (since each $Y_j$ is also finite or countable). For each $Y_j\in C_v$, define $\widetilde{Y_j}=Y_j\cup \{e_0^j\}$ if $Y_j\notin S$ (where $e_0^j$ is only a symbol) and define $\widetilde{Y_j}=Y_j$ if $Y_j\in S$.

Our next goal is to define $R_e$ for each $e\in s^{-1}(v)$.

Partition the interval $D_v$ into $|\widetilde{Y_1}|$ intervals closed on the left and open on the right, and call the intervals $I_{e_i^1}$ where $e_i^1\in \widetilde{Y_1}$. For each $e_i^1\in Y_1$ define $R_{e_i^1}=I_{e_i^1}$.
Note that the set of intervals $\left\{I_{e_i^1}:e_i^1\in \widetilde{Y_1}\right\}$ is a countable set, and so we may write it as $\{^1I^k:k\in \N\}$.

Now, partition each interval $^1I^k$ into $|\widetilde{Y_2}|$ closed on the left and open on the right intervals, called $^1I_{e_j^2}^k$, where $e_j^2\in \widetilde{Y_2}$. Define, for each $e_j^2\in Y_2$, $$R_{e_j^2}=\bigcup\limits_{k\in \N}{}^1I_{e_j^2}^k.$$

Since $\left\{^1I_{e_j^2}^k:k\in \N,\,\,\, e_j^2\in \widetilde{Y_2}\right\}$ is countable, we may write this set as $\left\{ ^2I^k:k\in \N\right\}$.

Partition each interval $ ^2I^k$ into $|\widetilde{Y_3}|$ (open on the right and closed on the left) intervals $^2I_{e_j^3}^k$. Define, for each $e_j^3\in Y^3$, $$R_{e_j^3}=\bigcup\limits_k {}^2I_{e_j^3}^k.$$

In general, given a partition $\{^nI^k:k\in \N\}$ of $D_v$ (obtained as above), partition each interval $^nI^k$ in $|\widetilde{Y_{n+1}}|$ closed on the left and open on the right intervals $^nI_{e_j^{n+1}}^k$ where $e_j^{n+1}\in \widetilde{Y_{n+1}}$. Then, for each $e_j^{n+1}\in Y_{n+1}$, define $$R_{e_j^{n+1}}=\bigcup\limits_{k}{}^nI_{e_j^{n+1}}^k.$$

So, we obtain $R_e$, for each $e\in s^{-1}(v)$. By applying this process to each non-sink $v\in E^0$ we obtain $R_e$ for all $e\in E^1$.
It is not hard to see that the sets $R_e$ satisfy the conditions $1,2$ and $3$ of the theorem.

To obtain the desired branching system, define $\X=\bigcup\limits_{v\in E^0}D_v$. It is also not hard to see that the families $\{R_e\}_{e\in E^1}$, $\{D_v\}_{v\in E^0}$ satisfy the relations 1-4 from Definition \ref{brancsystem}.
Finally, we need to obtain bijections $f_e:D_{r(e)}\rightarrow R_e$ for all $e\in E^1$. Fix $e\in E^1$.  By the definition of $R_e$, we see that $R_e$ is a union of closed on the left and open on the right disjoint intervals $\{J_k\}_{k\in \Delta}$, where $\Delta$ is finite or countable. Partition $D_{r(e)}$ into $|\Delta|$ closed on the left and open on the right (disjoint) intervals $D_k$ with $k\in \Delta$. Then since, for each $k\in \Delta$, let $f_k:D_k\rightarrow J_k$ be a bijective map (for example, the linear map). Now, given $x\in D_{r(e)}$, then $x\in D_k$ for some $k\in \Delta$, and define $f_e(x):=f_k(x)$. Then $f_e:D_{r(e)}\rightarrow R_e$ is a bijective map.
\fim

\begin{obs}\label{r1} It is not hard to see in the previous proof that if $X_1,...,X_n,...,X_m\in C_v$ are disjoint sets with $X_{n+1},...,X_m\notin S$ and if $e_i\in X_i$ for $1\leq i\leq n$ then

$$R_{e_1}\cap...\cap R_{e_n}\cap (D_v\setminus\bigcup\limits_{e\in X_{n+1}}R_e)\cap...\cap (D_v\setminus\bigcup\limits_{e\in X_m}R_e)\neq \emptyset$$

\end{obs}

\begin{exemplo}\label{ex1} The graph of this example is a graph with $4$ edges, without loops, and with injective range, as follows:

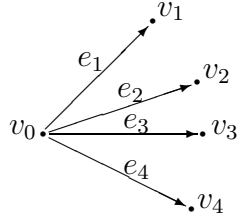
\centerline{
\setlength{\unitlength}{1,5cm}
\begin{picture}(2,1.5)
\put(0,0){\circle*{0.05}}
\put(-0.3,0){$v_0$}
\put(1.4,0){\circle*{0.05}}
\put(0.05,0){\vector(1,0){1.3}}
\put(1.35,0.46){\circle*{0.05}}
\put(0.05,0.02){\vector(3,1){1.25}}
\put(0.96,1){\circle*{0.05}}
\put(0.02,0.05){\vector(1,1){0.9}}
\put(0.05,-0.03){\vector(2,-1){1.2}}
\put(1.3,-0.66){\circle*{0.05}}
\put(1.35,-0.65){$v_4$}
\put(0.7,0.055){$e_3$}
\put(0.7,-0.35){$e_4$}
\put(0.65,0.33){$e_2$}
\put(0.3,0.6){$e_1$}
\put(1,1.05){$v_1$}
\put(1.4,0.5){$v_2$}
\put(1.47,0){$v_3$}
\end{picture}}
\end{exemplo}
\vspace{0.8cm}

Let $X_1=\{e_1,e_2\}$ and $X_2=\{e_3,e_4\}$, and let $S=\{X_2\}$. We follow the proof of the previous theorem to obtain $D_{v_i}$ and $R_{e_j}$. Set $D_{v_i}=[i,i+1)$ for $0\leq i\leq 5$.
To obtain $R_{e_j}$, proceed as follows:
\begin{itemize}
\item since $X_1\notin S$, partition the interval $[0,1)$ into 3 intervals, \newline\centerline{$[0,1)=[0,\frac{1}{3})\cup[\frac{1}{3},\frac{2}{3})\cup[\frac{2}{3},1).$} Define $R_{e_1}=[0,\frac{1}{3})$ and $R_{e_2}=[\frac{1}{3},\frac{2}{3})$

\item since $X_2 \in S$, partition the 3 intervals into 2 intervals, as follows:
\newline\centerline{$[0,\frac{1}{3})=[0,\frac{1}{6})\cup [\frac{1}{6},\frac{1}{3}),$}
\newline\centerline{$[\frac{1}{3},\frac{2}{3})=[\frac{1}{3},\frac{1}{2})\cup [\frac{1}{2},\frac{2}{3}),$}
\newline\centerline{$[\frac{2}{3},1)=[\frac{2}{3},\frac{5}{6})\cup [\frac{5}{6},1).$}
Define
\newline\centerline{$R_{e_3}=[0,\frac{1}{6})\cup [\frac{1}{3},\frac{1}{2})\cup [\frac{2}{3},\frac{5}{6})$} and \newline\centerline{$R_{e_4}=[\frac{1}{6},\frac{1}{3})\cup [\frac{1}{2},\frac{2}{3})\cup [\frac{5}{6},1).$}
\end{itemize}

By Remark \ref{r1}, since $X_1\notin S$,  for $h\in X_2$ it holds that $R_h\cap (D_{v_0}\setminus\bigcup\limits_{e\in X_1} R_e)\neq \emptyset.$ For example (if $h=e_3$) $R_{e_3}\cap (D_{v_0}\setminus\bigcup\limits_{e\in X_1} R_e)=[\frac{2}{3},\frac{5}{6})$.

Theorem \ref{existencebrancsys} together with theorem \ref{rep} guarantees that every Cohn-Leavitt path algebra of separated graphs $L_K(E,C,S)$ of a countable graph $E$  may be represented in $Hom_K(M)$. Let us summarize this result in the following corollary:

\begin{corolario}\label{cor1} Given a triple $(E,C,S)$, with $E$ countable, there exists a homomorphism $\pi:L_K(E,C,S)\rightarrow Hom_K(M)$ such that $$\pi(v)(\phi)=\chi_{D_v}.\phi,\,\,\,\,\,\pi(e)(\phi)=\chi_{R_e}.\phi\circ f_e^{-1}\,\,\,\,\text{ and }\,\,\,\,\pi(e^*)(\phi)=\chi_{D_{r(e)}}.\phi\circ f_e$$ for each $\phi\in M$,  where $M$ is the $K$ module of all functions from $\X$ taking values in $K$, $\X$ is an (possible unlimited) interval of $\R$, and $R_e$ and $D_v$ are as in theorem \ref{existencebrancsys}
\end{corolario}

\begin{corolario} In the algebra $L_K(E,C,S)$ it holds that:

\begin{enumerate}
\item $e\neq 0$ for each $e\in E^1$,
\item $v\neq 0$ for each $v\in E^0$,
\item $e^*f\neq 0$ for each $e\in X$, $f\in Y$, $X,Y\in C_v$ with $X\neq Y$,
\item for each finite set $X\in C_v\setminus S$, $$\sum\limits_{e\in X}ee^*v=\sum\limits_{e\in X}ee^*=v\sum\limits_{e\in X}ee^*$$ but $\sum\limits_{e\in X}ee^*\neq v$.
\item for each finite sets $X,Y\in C\setminus S$ it holds that $\sum\limits_{e\in X}ee^*\neq \sum\limits_{f\in Y}ff^*$.
\end{enumerate}

\end{corolario}

\demo Consider the homomorphism $\pi:L_K(E,C,S)\rightarrow Hom_K(M)$ as in the previous corollary. Since $\pi(v)(\phi)=\chi_{D_v}\cdot \phi$ for each $\phi \in M$ then $\pi(v)\neq 0$, for each $v\in E^0$ and so $v\neq 0$ in $L_K(E,C,S)$. Moreover, $\pi(e)\pi(e^*)(\phi)=X_{R_e}\cdot \phi$, and so $\pi(e)\neq 0$ and also $e\neq 0$ in $L_K(E,C,S)$.

Item 3 follows by item 1 of the previous theorem. In fact, note that $\pi(e)\pi(e^*)\pi(f)\pi(f^*)(\phi)=\chi_{R_e\cap R_f}\cdot \phi$, and since $R_e\cap R_f\neq \emptyset$ (for $e\in X$, $f\in Y$, $X,Y\in C_v$ and $X\neq Y$) then $\pi(e)\pi(e^*)\pi(f)\pi(f^*)\neq 0$, and so $e^*f\neq 0$.

Let us prove item 4. The equalities $$\sum\limits_{e\in X}ee^*v=\sum\limits_{e\in X}ee^*=v\sum\limits_{e\in X}ee^*$$ follow by $(E1)$ and $(E2)$ from the definition of $L_K(E,C,S)$, and $\sum\limits_{e\in X}ee^*\neq v$ follows from the second item of the previous theorem.

To prove item 5, let $X,Y\in C$ be finite sets with $X\neq Y$. If $X\in C_u$ and $Y\in C_v$ ($u\neq v$) then $$\sum\limits_{e\in X}ee^*\sum\limits_{f\in Y}ff^*=\sum\limits_{e\in X}ee^*uv\sum\limits_{f\in Y}ff^*=0,$$ and since $\sum\limits_{e\in X}ee^*\neq 0$ then $\sum\limits_{e\in X}ee^*\neq\sum\limits_{f\in Y}ff^*$. If $X,Y\in C_v$ then the inequality $\sum\limits_{e\in X}ee^*\neq\sum\limits_{f\in Y}ff^*$ follows from the third item of the  previous theorem.
\fim

\section{Injectivety and the amalgamated free product structure}

The representations introduced in the previous sections are  adaptations of the representations introduced and studied in \cite{rg1}, \cite{perroncuntz}, \cite{repgraph} and \cite{uniteq} to the separated graph case. But, contrary to what happened to Leavitt algebras, for most graphs, representations arising from branching systems on separated graphs can not be faithful, since for edges $e$ and $f$ in different sets of a partition of a vertex $v$, we always have that $S_eS_e^*$ commute with $S_fS_f^*$, but this in general is not true in the algebra. In order to present some examples of injective representations we have then to look at the abelianized algebra, which is a quotient of $L_K(E,C,S)$.

\begin{definicao}(as in \cite{AE}) Let $(E,C)$ be a separated graph. The abelianized Cohn-Leavitt algebra, denoted by $AL_K(E,C,S)$ is the quotient of the Cohn-Leavitt algebra $L_K(E,C,S)$ by the ideal $J$ generated by all the elements $\lambda\lambda^*\beta\beta^*-\beta\beta^*\lambda\lambda^*$, where $\lambda,\beta$ belong to the multiplicative semigroup generated by $E^1\cup (E^1)^*$.

%universal $K$-algebra generated by a set
%$\{v:v\in E^0\}$, of pairwise orthogonal idempotents, together with a set $\{e,e^*:e\in E^1\}$ of elements satisfying the relations E1, E2, SCK1, SCK2 of the definition of $L_K(E,C,S)$ and also satisfying:

%\begin{itemize}
%\item[(A1)] $\lambda \lambda^* \beta \beta^* = \beta \beta^* \lambda \lambda^* $ for all paths $\lambda$, $\beta$ in $\text{Path}(E)$.
%\end{itemize}
\end{definicao}

{\bf Remark:} Notice that any representation $\pi$ of $L_K(E,C,S)$ arising from a branching system, as in Theorem \ref{rep}, is automatically a representation of $AL_K(E,C,S)$, since for each $\beta,\gamma$ in the multiplicative semigroup generated by $E^1\cup (E^1)^*$, $\pi(\gamma\gamma^*)$ and $\pi(\beta\beta^*)$ are multiplication operators in $Hom_K(M)$, and so $\pi(\gamma\gamma^*)\pi(\beta\beta^*)-\pi(\beta\beta^*)\pi(\gamma\gamma^*)=0$.

\vspace{1pc}
Next we show that, for a class of graphs, any representation of $AL_K(E,C,S)$ arising from a branching system as in theorem \ref{existencebrancsys} is faithful. More specifically we will consider countable graphs with no loops where all edges have the same source and the range map is injective. An example of such a graph was given in example \ref{ex1}.

In order to proceed we need to recall the description of a basis for $L_K(E,C,S)$. This was done in \cite{ara1} by P. Ara and K. R. Goodearl. For the graphs in question, a basis for $L_K(E,C,S)$ consists of the set $\mathcal{B}$ of paths of the form 

$$\alpha:= \mu e_1^*e_2e_2^*\ldots e_{n-1}e_{n-1}^*e_n \nu^*,$$

such that $\alpha$ is $C$ separated and reduced with respect to $S$ (and $\mu$ and $\nu$ are allowed to have length zero). Notice that for the case in mind the source of all edges is a vertex $v$, and so $\alpha$ is $C$-separated iff $e_i$ and $e_{i+1}$ are in different sets $X,Y \in C_v$ for all $i$. Furthermore, $\alpha$ is reduced with respect to $S$ iff for each $X \in S$, an edge $e_X$ has been selected and $e_ie_i^* \neq e_Xe_X^*$, for all $i$, or if $\alpha$ is equal to an edge $e$, a ghost edge $e^*$ or a vertex.

We can now prove the faithfulness of our representations, but firs we need the following lemma.

\begin{lema}\label{lemma1} Let $(E,C,S)$ be a separated graph, where all edges have the same source, $v_0$, and the range map is injective. Let $\pi$ be the representation of $AL_K(E,C,S)$ arising from the branching system defined in theorem \ref{existencebrancsys}. Let $x\in AL_K(E,C,S)$ be a non zero linear combination of elements of the form $x=\gamma_0v_0 + \sum_{j=1}^N \gamma_j e_1^j(e_1^j)^*\ldots e_{n_j}^j(e_{n_j}^j)^*$, where $e_i^j$ are edges. Then $\pi(x) \neq 0$.
\end{lema}

\demo First note that we may suppose, for each $j$, that $e_i^j$ and $e_k^j$ are $C$-separated for each $i\neq k$ (otherwise $(e_i^j)^*e_k^j=0$). Moreover, we may suppose for each $j$ that $e_1^j(e_1^j)^*\ldots e_{n_j}^j(e_{n_j}^j)^*$
is reduced with respect to $S$. \footnote{that is, for all $i,j$ and for each $X\in S$, we may suppose that $e_i^j\neq e_X$ (where $e_X\in X$ has been previously selected to obtain $B$) by replacing $e_i^j(e_i^j)^*$ by $v_0-\sum\limits_{e\in X\setminus\{ e_X\}}ee^*$ if $e_i^j=e_X$ for some $i,j$ and $X\in S$.}

If $\gamma_j=0$ for each $j$ then $\pi(x)=\gamma_0\pi(v_0)=\gamma_0(1_{v_0})\neq 0$. So, suppose $\gamma_j\neq 0$ for each $1\leq j\leq N$.

Let $X_1,...,X_m\in C_{v_0}$ be the subsets which contain some edge $e_i^j$, for $1\leq j\leq N$ and $1\leq i\leq n_j$.

Case 1: Suppose $\gamma_0\neq 0$.

For each $1\leq k\leq m$, define $y_k$ as follows:
if $X_k$ is infinite, let $y_k=e_ke_k^*$ where $e_k\neq e_i^j$ for each $1\leq j\leq N$ and $1\leq i\leq n_j$; if $X_k$ is finite and $X\in S$ leq $y_k=e_{X_k}e_{X_k}^*$ (where $e_{X_k}e_{X_k}^*$ has been selected to form $\mathcal{B}$); and if $X_k$ is finite and $X\notin S$ let $y_k=v_0-\sum\limits_{e\in X_k}ee^*$. So, for each $e_i^j\in X_k$ it holds that $e_i^j(e_i^j)^*y_k=0$. In particular, for each $1\leq j\leq N$ and $1\leq i\leq n_j$,

$$e_i^j(e_i^j)^*y_1...y_m=0,$$ and therefore, for each $1\leq j\leq N$,

$$e_1^j(e_1^j)^*...e_{n_j}^j(e_{n_j}^j)^*y_1...y_m=0.$$

So, by multiplying the equality
$$x=\gamma_0v_0 + \sum_{j=1}^N \gamma_j e_1^j(e_1^j)^*\ldots e_{n_j}^j(e_{n_j}^j)^*$$
 by $y_1...y_m$ we obtain

$$xy_1...y_m=\gamma_0v_0y_1...y_m=\gamma_0y_1...y_m.$$

and then

$$\pi(xy_1...y_m)=\gamma_0\pi(y_1....y_m)=\gamma_0\pi(y_1)...\pi(y_m).$$ Since $\pi(y_k)=1_{R_{e_k}}$ for some $e_k\in X_k$ or $\pi(y_k)=1_{D_{v_0}\setminus\bigcup\limits_{e\in X_k}R_e}$, then, by Remark \ref{r1},  $\pi(y_1)...\pi(y_m)\neq 0$. So, we have $\pi(x)\pi(y_1...\pi_m)=\pi(xy_1...y_m)=\gamma_0\pi(y_1...y_m)\neq 0$, and then $\pi(x)\neq 0$.

Case 2: Suppose $\gamma_0=0$.

Choose some $t$, $1\leq t\leq N$, such that $n_t\leq n_j$ for all $1\leq j\leq N$. For each $1\leq k\leq m$  define $y_k$ as follows: if $e_i^t\in X_k$ for some $1\leq i\leq n_t$ define $y_k=e_i^t(e_i^t)^*$, and otherwise define $y_k$ as in Case 1. So, for each $j\neq t$, there is an $1\leq i\leq n_j$ such that $e_i^j\neq e_l^t$ for all $1\leq l\leq n_t$. Let $k$ be such that $e_i^j\in X_k$. Then $e_i^j(e_i^j)^*y_k=0$. As in Case 1,
$$e_1^j(e_1^j)^*...e_{n_j}^j(e_{n_j}^j)^*y_1...y_m=0$$

for each $1\leq j\leq N$ with $j\neq t$, and then

$$xy_1...y_m=\gamma_te_1^t(e_1^t)^*...e_{n_t}^t(e_{n_t}^t)^*y_1...y_m.$$

Since for each $1\leq i\leq n_t$ the element $e_i^t(e_i^t)^*$ equals to some $y_k$ then $$e_1^t(e_1^t)^*...e_{n_t}^t(e_{n_t}^t)^*y_1...y_m=y_1...y_m,$$ and so
$$xy_1...y_m=\gamma_ty_1...y_m.$$

Then it follows, as in Case 1, that $\pi(x)\neq 0$.
\fim

\begin{teorema}\label{teorema1}
Let $(E,C,S)$ be a separated graph, where all edges have the same source, $v_0$, the range map is injective and $E$ has no loops. Then the representation $\pi$ of $AL_K(E,C,S)$ arising from the branching system in theorem \ref{existencebrancsys} is faithful.
\end{teorema}

\demo
Let $x\in AL_K(E,C,S)$ be a nonzero element. By \cite[2.7]{ara1}, 
$$x=\gamma_0v_0+\sum_{j=1}^k \gamma_j \alpha_j $$ where $a_j\in \mathcal{B}$ for each $j$, and since the unique paths in $E$ are the edges, then $\alpha_j=\mu_j(e_1^j)^*e_2^j(e_2^j)^*\ldots e_{n_j-1}^j(e_{n_j-1}^j)^*e_{n_j}^j\nu_j^*$ for each $j$, where each $e_i^j$ is an edge and each $\nu_j$ and each $\mu_j$ has length zero or is an edge. We will show that $\pi(x)\neq 0$.

First notice that, since the set of finite sums of vertices is a set of local units for $AL_K(E,C,S)$, there exists vertices $v$ and $w$ such that $vxw\neq 0$. 

Suppose $v\neq v_0$ and $w=v_0$. By hypothesis, there exists only one edge $e$ such that $r(e)=v$. Writing $v=e^*e$ we obtain that $v\alpha_j\neq 0$ only if $|\mu_j|=0$ and $e_1^j=e$, and since $w=v_0$ then $\alpha_jw\neq 0$ only if $\nu_j=e_{n_j}^j$. Then 

$$0\neq vxw=\sum_{j:v\alpha_jw\neq 0} \gamma_je^*e_2^j(e_2^j)^*\ldots e_{n_j}^j(e_{n_j}^j)^*=$$
$$=e^*\sum_{j:v\alpha_jw\neq 0} \gamma_j ee^*e_2^j(e_2^j)^*\ldots e_{n_j}^j(e_{n_j}^j)^*.$$ By the previous lemma, 

$$\pi\left(\sum_{j:v\alpha_jw\neq 0} \gamma_j ee^*e_2^j(e_2^j)^*\ldots e_{n_j}^j(e_{n_j}^j)^*\right)\neq 0$$ and so 

$$0\neq \pi\left(\sum_{j:v\alpha_jw\neq 0} \gamma_j ee^*e_2^j(e_2^j)^*\ldots e_{n_j}^j(e_{n_j}^j)^*\right)=\pi(euxw)=\pi(e)\pi(u)\pi(x)\pi(w),$$ and so $\pi(x)\neq 0$.

The cases $v\neq v_0$ and $w\neq v_0$, $v=v_0$ and $w=v_0$, $v=v_0$ and $w\neq v_0$ 
follow in a similar way and are left to the reader.
\fim

We finalize the paper showing that the Cohn-Leavitt path algebra of a separated graph can be written as an amalgamated free product, with a common subset of idempotents, of Conh-Leavitt path algebras over graphs where all edges have the same source.

So, let $(E,C,S)$ be a triple as in Definition \ref{defleviatt}. For each $X\in C$ consider the directed graph $E_X=(E_X^0, X, r, s)$, where $E_X^0$ is a copy of $E^0$, that is, $E_X^0=\{v_X:v\in E^0\}$, and abusing the notation, $r,s$ are the range ans source maps of $E$ restricted to $X$. Notice that all edges in $E_X$ have the same source. Define, for each $X\in C$, $A_X$ as the universal $K$-algebra generated by $\{e,e^*:e\in X\}\cup E_X^0$ with relations given by:
\begin{enumerate}
\item the elements of $E_X^0$ are pairwise orthogonal idempotents,
\item $er(e)=e=s(e)e$ and $r(e)e^*=e^*=e^*s(e)$, for all $e\in X$,
\item $e^*f=\delta_{e,f}r(e)$, for all $e,f\in X$,
\item $v_X=\sum\limits_{e\in X}ee^*$, if $X\in S$ and $v_X=s^{-1}(X)$.
\end{enumerate}

Notice that if $X\in S$, or if $X$ is infinite, then $A_X$ is the Leavitt path algebra of the directed graph $E_X$. In particular, if $S=C_{fin}$ then each $A_X$ is an Leavitt path algebra.

Now let $A$ be the free product of the $K$-algebras $A_X$, and let $I$ be the two sided ideal of $A$ generated by the set $\{v_X-v_Y:X,Y\in C; v\in E^0\}$. The quotient algebra $A/I$ is called the amalgamated free product of $\{A_X\}_{X\in C}$ with common subset $E^0$.

With the above in mind, we obtain the following proposition:

\begin{proposicao} The Cohn-Leavitt path algebra $L_K(E,C,S)$ of the triple $(E,C,S)$ is $K$-isomorphic to the amalgamated free product $A/I$.
\end{proposicao}

\demo The proof follows by using the universal property of  $L_K(E,C,S)$ to define a $K-$homomorphism $\psi:L_K(E,C,S)\rightarrow A/I$ such that $\psi(e)=[e]$, for each $e\in E^1$, and $\psi(v)=[v_X]$ (where $X$ is some set $X\in C$ and $[v_X]$ denotes the equivalence class of $v_X$), and then use the universal property of $A$ to define the inverse of $\psi$.

\fim

\begin{obs} The above result is an extension to the separated graph case of a result proved by Larki in the context of Leavitt path algebras of edge colored graphs (see \cite{larki}). As it happens, each edge-colored graph $G=(V,E,r,s,d)$ (in the sense of \cite{larki}) can be seen as a separated graph, with partitions $C_v=\{s^{-1}(v)\cap d^{-1}(i):i\in \N \text{ and } s^{-1}(v)\cap d^{-1}(i)\neq \emptyset\}$, and so the Leavitt path algebra of the edge colored graph $G$ (as defined in \cite{larki}) and the Leavitt path algebra of the above separated graph coincide.
\end{obs}

\begin{corolario} Let $(E,C,S)$ be a separated graph with no loops and such that the range map is injective. Then $AL_K(E,C,S)$ can be written as an amalgamated free product of abelizanized Cohn-Leavitt path algebras that can be faithfully represented as in theorem \ref{teorema1}.

\end{corolario}

\vspace{1.5pc}

Daniel Gonçalves, Departamento de Matemática, Universidade Federal de Santa Catarina, Florianópolis, 88040-900, Brasil

Email: daemig@gmail.com

\vspace{0.5pc}
Danilo Royer, Departamento de Matemática, Universidade Federal de Santa Catarina, Florianópolis, 88040-900, Brasil

Email: danilo.royer@ufsc.br
\vspace{0.5pc}

\end{document}